\documentclass[12pt]{article}
\usepackage{amsmath,amsthm,amssymb}
\usepackage{mathrsfs}
\usepackage[titletoc]{appendix}
\usepackage{bbm}
\usepackage{float}
\usepackage{tikz}
\usepackage[all]{xypic}
\usepackage{graphicx}
\usepackage{latexsym,amsfonts,amssymb,amsmath,amscd}
\usepackage{amsfonts}
\newcommand{\End}{{\rm End}\,}
\begin{document}
\input amssym.def
\newcommand{\singlespace}{
    \renewcommand{\baselinestretch}{1}
\large\normalsize}
\newcommand{\doublespace}{
   \renewcommand{\baselinestretch}{1.2}
   \large\normalsize}
\renewcommand{\theequation}{\thesection.\arabic{equation}}

\setcounter{equation}{0}
\def \ten#1{_{{}_{\scriptstyle#1}}}
\def \Z{\Bbb Z}
\def \C{\Bbb C}
\def \R{\Bbb R}
\def \Q{\Bbb Q}
\def \N{\Bbb N}
\def \F{\Bbb F}
\def \L{\Bbb L}
\def \l{\lambda}
\def \V{V^{\natural}}
\def \wt{{\rm wt}}
\def \Spec{{\rm Spec}}
\def \Span{{\rm Span}}
\def \tr{{\rm tr}}
\def \Res{{\rm Res}}
\def \End{{\rm End}}
\def \Im{{\rm Im}}
\def \Aut{{\rm Aut}}
\def \mod{{\rm mod}}
\def \Id{{\rm Id}}
\def \Hom{{\rm Hom}}
\def \im{{\rm im}}
\def \<{\langle}
\def \>{\rangle}
\def \w{\omega}
\def \c{{\tilde{c}}}
\def \o{\omega}
\def \t{\tau }
\def \ch{{\rm ch}}
\def \a{\alpha }
\def \b{\beta}
\def \e{\epsilon }
\def \la{\lambda }
\def \om{\omega }
\def \O{\Omega}
\def \qed{\mbox{ $\square$}}
\def \pf{\noindent {\bf Proof: \,}}
\def \voa{vertex operator algebra\ }
\def \voas{vertex operator algebras\ }
\def \p{\partial}
\def \1{{\bf 1}}
\def \ll{{\tilde{\lambda}}}
\def \H{{\bf H}}
\def \F{{\Bbb F}}
\def \h{{\frak h}}
\def \g{{\frak g}}
\def \rank{{\rm rank}}
\def \({{\rm (}}
\def \){{\rm )}}
\def \Y {\mathcal{Y}}
\def \I {\mathcal{I}}
\def \A {\mathcal{A}}
\def \B {\mathcal {B}}
\def \Cc {\mathcal {C}}
\def \H {\mathcal{H}}
\def \M {\mathcal{M}}
\def \V {\mathcal{V}}
\def \O{{\bf O}}
\def \AA{{\bf A}}
\def \1{{\bf 1}}
\def\Ve{V^{0}}
\def\ha{\frac{1}{2}}
\def\se{\frac{1}{16}}
\def\Der{{\rm Der}}
\def\qdim{{\rm qdim}}
\singlespace
\newtheorem{Theorem}{Theorem}[section]
\newtheorem{Proposition}[Theorem]{Proposition}
\newtheorem{Lemma}[Theorem]{Lemma}
\newtheorem{Corollary}[Theorem]{Corollary}
\newtheorem{Remark}[Theorem]{Remark}
\newtheorem{Conjecture}[Theorem]{Conjecture}
\newtheorem*{CPM}{Theorem}
\newtheorem{Definition}[Theorem]{Definition}

\begin{center}
{\Large {\bf Hopf actions on vertex operator algebras}} \\

\vspace{0.5cm} Chongying Dong\footnote{Supported by a NSF grant}
\\
 Department of Mathematics, University of
California, Santa Cruz, CA 95064 USA \\
Hao Wang\\
 School of Mathematics,  Sichuan University,
Chengdu 610064 China
\end{center}
\hspace{1.5 cm}

\begin{abstract}
The Hopf actions on vertex operator algebras are investigated.   If the action is semisimple, a  Schur-Weyl type decomposition is obtained. When the Hopf algebra is finite dimensional and the action is faithful,   the action is a group action.  Moreover  if the Hopf algebra is finite dimensional and the action is semisimple and inner faithful, the action is also a group action.  In this case, inner faithfulness is equivalent to faithfulness.
\end{abstract}

\section{Introduction}

The notion of Hopf action on a vertex operator algebra studied in this paper unifies and generalizes  the notions of automorphisms and
derivations of a vertex operator algebra. Although the generalization is natural,  our study of Hopf actions on vertex operator algebras
is motivated by the Schur-Weyl type duality for the action of a compact Lie group  on a vertex operator algebra \cite{DLM3, DM}
and the inner faithful action of a finite dimensional semisimple Hopf algebra on a commutative domain \cite{EW1}. 

Consideration of a group action $G$ on a vertex operator algebra $V$ is a central problem in the theory of vertex operator algebras. The well known orbifold theory conjecture says if $V$ is rational (the admissible or $\Z_+$-graded module category is semisimple \cite{DLM2, LL}) and $G$ is finite then the fixed point vertex operator algebra $V^G$ is also rational and each irreducible $V^G$-module occurs in an irreducible $g$-twisted $V$-module (See \cite{FLM, DLM2} for the definition of twisted module). The rationality of $V^G$ has been established recently in \cite{CM} if $G$ is solvable and $V$ is $C_2$-cofinite and the second part of the conjecture  has been proved  with the assumption that $V$ is $C_2$-cofinite and $V^G$ is rational in \cite{DRX}. So decomposition of an irreducible $g$-twisted module as direct sum of irreducible $V^G$-modules becomes important in understanding the $V^G$-module category. This was first achieved for $V$-module $V$ for any compact Lie group $G$ in \cite{DLM3}. A Schur-Weyl type duality 
was obtained as follows:
$$V=\oplus_{i\in I}M_i\otimes V_{M_i}$$
 where $\{M_i|i\in I\}$ is the set of inequivalent finite dimensional irreducible  $G$-modules and $V_{M_i}$ is the multiplicity space of
$M_i$ in $V.$ Moreover, each  $V_{M_i}$ is an irreducible  $V^G$-module and $V_{M_i},V_{M_j}$ are isomorphic if and only if $i=j$.
That is, $(G,V^G)$ forms a dual pair in the sense of \cite{H}.  A decomposition of an irreducible $g$-twisted $V$-module into direct sum of irreducible $V^G$-modules is also carried out in \cite{DRX}. The first motivation of our work is to give a similar Schur-Weyl type duality for the Hopf action on a simple vertex operator algebra.

The actions of Hopf algebras on rings and various algebras have a long history \cite{CWWZ, CWZ, EW1, EW2, M}.  It is proved in \cite{EW1} that any finite dimensional semisimple and inner faithful action  on a commutative domain  is a group action. A commutative domain $A$ has two  properties: (1)  If  $ 0\neq a,0\neq b\in A$ then  $ab\neq0,$ (2) $[a,b]=ab-ba=0$ for $ a,b\in A.$ Although a vertex operator algebra is not a commutative domain,  a  simple VOA $V$  has the  similar properties\cite{DL}: 

   (1) If  $0\neq u,0\neq v\in V,$ then $Y(u,z)v\neq0,$

 (2) If $u,v\in V$ then there exists $n\in \mathbb{Z}_+$ such that  $(z_1-z_2)^n[Y(u,z_1),Y(v,z_2)]=0.$\\
So it is natural to expect that the result in \cite{EW1} holds for the semisimple and inner faithful Hopf action on a simple vertex operator algebra. This is our second motivation.
 
We now discuss our main results and the ideas of the proofs.  First, let $H$ be a Hopf algebra and $V$ a simple vertex operator algebra such that $V$ is a semisimple $H$-module. Similar to the group action case \cite{DLM3}, we have 
a Schur-Weyl type decomposition of $V$ as a module for the pair $(H, V^H):$
 $$V=\oplus_{i\in I}M_{i}\otimes V_{M_i}$$
 where $\{M_i|i\in I\}$ is the set of inequivalent finite dimensional irreducible representations of $H$ occurring in $V$ and 
the multiplicity space $ V_{M_i}$ is an irreducible $V^H$-module. For different $i,j$, $V^H$-modules  $V_{M_i}$ and $V_{M_j}$ are inequivalent $V^H$-modules. We should point out that in the group action case in \cite{DLM3}, every finite dimensional irreducible $G$-module occurs, but this is not true anymore in the current situation. For example, if $H$ is the universal enveloping algebra of a simple Lie algebra $\frak{g}$ and $V$ is the affine vertex operator algebra, not every irreducible $\frak{g}$-module occurs. The main tool used to prove this result is the $A_n(V)$ (which is an associative algebra)  theory developed in \cite{DLM2, Z}.  Let $V_{M_i}=\oplus_{n\geq 0}(V_{M_i})_n$ where
$(V_{M_i})_n$ is the subspace of $V_{M_i}$ consisting of vectors of weight $n.$ According to \cite{DLM2}, $V_{M_i}$ is an irreducible
$V^H$-module if and only if $(V_{M_i})_n$ is an irreducible $A_n(V^H)$-module for every $n.$ We established the irreducibility of
$(V_{M_i})_n$ as an  $A_n(V^H)$-module. 

 Second, we consider the case when $H$ is a finite dimensional Hopf algebra and the action is faithful. Using the properties of vertex 
operator algebra to show  that $\Delta(h)Y(u,z)v=\Delta^{op}(h)Y(u,z)v$ for $h\in H$ and $u,v\in V$, and $H$ is cocommutative. Thus $H$ is a group algebra. 

Finally, we consider the inner faithful action. In this case we assume $H$ is finite dimensional and the action is semisimple and inner faithful. The main idea is to prove the kernel $K$ of the action is a Hopf ideal by using the following properties for a simple vertex operator algebra:

(i)  $V$ is spanned by $\{u_nv|u\in V,n\in \mathbb{Z}\}$ for any fixed $0\neq v\in V$  \cite{DM,L},

 (ii) If $u^1,\cdots,u^s\in V$ are not all zero, and $v^1,\cdots,v^s\in V$ are linearly independent. then, $\sum_{i=1}^sY(u^i,z)v^i\neq0$ \cite{DM}. 
Inner faithful assumption implies that $K=0.$ So the action of $H$ is faithful and a group action. Although we believe that the kernel of
a Hopf action on a vertex operator algebra is a Hopf ideal but we could not prove it in this paper. 
     
There is a related problem on the Hopf action. That is, classify the irreducible $V^H$-modules. From the discussion above, we know that
if $H$ is a group algebra of a finite group, we need twisted modules. But for general $H$ we do not know what the analogue of the twisted module is. 
  
This paper is organized as follows:  In Section 2, we review some basic and well known results  about Hopf algebras.  In section 3, after reviewing  some facts about vertex operator algebras and modules, we  define the Hopf actions on vertex operator algebras, that is, $H$-module vertex operator algebras for Hopf  algebra $H.$ Then we establish the Schur-Weyl type decomposition of a simple vertex operator algebra $V$ as $(H,V^H)$-module if $V$ is a semisimple $H$-module. In section 4, we discuss the fusion rules among
irreducible $V^H$-modules occurring in $V$ in terms of the tensor product of the corresponding  simple $H$-modules.
In section 5, we consider the case when $H$ is finite dimensional and the Hopf action is faithful, and  prove that $H$ is a group algebra. In section 6, we remove the condition of faithfulness and assume $V$ is a semisimple $H$-module. Let $K$ be the kernel of the Hopf action, then we prove $K$ is a Hopf ideal. This implies that inner faithfulness is equivalent to faithfulness. In this case we show again that
the Hopf action is a group action. 

  We assume that the reader is familiar with the basics on VOA as presented in \cite{DLM,DLM1,DLM2,DLM3,FHL,FLM,LL}. We also refer the reader to \cite{DLM2} for various notions of  modules for a VOA and related results, \cite{FLM} for the  definition and properties of intertwining operators for $V$-modules, and \cite{A,DNR,M,M1} for  basic and important concepts in the theory of  Hopf algebras.   In this paper, we work over complex field $\mathbb{C}.$

We thank Professor Siu-Hung Ng for valuble suggestions and discussions on the Hopf actions on vertex operator algebras. 

\section{Basics on Hopf algebras}

  In this section, we recall from \cite{A,DNR,M} the basic facts on Hopf algebras.  Let $H=(H,\mu,\eta,\Delta,\epsilon,S)$ be a Hopf algebra where the linear maps $\mu:H\otimes H\rightarrow H,~\eta:\mathbb{C}\rightarrow H,~\Delta:H\rightarrow H\otimes H,~\epsilon:H\rightarrow \mathbb{C},~S:H\rightarrow H$ are multiplication, unit, comultiplication, counit and antipode, respectively. Then

  (i) $(H,\mu,\eta)$ is an associative algebra, i.e. for any $g,h,f\in H,$
  \begin{eqnarray}
  &&\mu(\mu\otimes \Id)(g\otimes h\otimes f)=\mu(\Id\otimes\mu)(g\otimes h\otimes f)=ghf,\\
  &&\mu(\eta\otimes \Id)=\mu(\Id\otimes\eta)=\Id.
  \end{eqnarray}

  (ii) $(H,\Delta,\epsilon)$ is a coassociative coalgebra, i.e. for any $h\in H,$
  \begin{eqnarray}\label{whh}
  &&(\Delta\otimes \Id)\Delta(h)=(\Id\otimes\Delta)\Delta(h),\\
  &&(\Id\otimes\epsilon)\Delta=(\epsilon\otimes \Id)\Delta=\Id.\label{ww}
  \end{eqnarray}
  In this paper, we always use Sweedler's notation and write $\Delta(h)=\sum h_1\otimes h_2$, by equation (\ref{whh}), we have
  \begin{eqnarray*}
  &&(\Delta\otimes Id)\Delta(h)=\sum\sum(h_1)_1\otimes (h_1)_2\otimes h_2\\
  &=&(Id\otimes\Delta)\Delta(h)=\sum\sum h_1\otimes (h_2)_1\otimes (h_2)_2\\
  &=&\sum h_1\otimes h_2 \otimes h_3.
  \end{eqnarray*}

  (iii) $\Delta,\epsilon$ are homomorphisms of associative algebras.

  (iv) $\mu,\eta$ are homomorphisms of coassociative coalgebras.

  (v) $S$ is the convolution inverse of the identity map $\Id$ of $H$, i.e. $$\mu(S\otimes \Id)\Delta=\mu(\Id\otimes S)\Delta=\eta\epsilon.$$
  $S$ is called the antipode of $H$.

  \begin{Remark}\label{0000}
  (1) In the above definition, (iii) is equivalent to (iv). See \cite{A} for the details of the proofs.

  (2) The convolution inverse means that for any $h\in H$, $$\sum h_1S(h_2)=\sum S(h_1)h_2=\epsilon(h).$$

  (3) By equation (\ref{ww}), for any $h\in H$, $\sum h_1\epsilon(h_2)=\sum\epsilon(h_1)h_2=h.$ So, if we define $\Delta_1=\Delta,\Delta_n=((\Id)^{\otimes^{n-1}}\otimes\Delta)\Delta_{n-1}$, for $n>1$. Then, for any $ i\leq n$, we have $((\Id)^{\otimes^i}\otimes\epsilon\otimes(\Id)^{\otimes^{n-i}})\Delta_n=\Delta_{n-1}.$

  (4) $S$ is an anti-endomorphism of $H$ as associative algebra, i.e. $S(gh)=S(h)S(g)$, for any $g,h\in H$(See \cite{A}).

  (5) For any $h\in H$, $\epsilon(S(h))=\epsilon(h),\Delta(S(h))=\sum S(h_2)\otimes S(h_1)$(See \cite{A}).
  \end{Remark}

  Since we work over complex field $\mathbb{C}$, the following results are hold.
  \begin{Lemma}\cite{DNR}\label{A,DNR}
  Let $H$ be a finite dimensional Hopf algebra and $S$ be the antipode of $H$. Then, $S$ is of finite order, i.e. there $n\in \N$, such that $S^n=\Id$.
  \end{Lemma}

  \begin{Lemma}\cite{DNR}
  Let $H$ be a finite dimensional Hopf algebra, if $H$ is cocommutative, then $H$ is a group algebra.
  \end{Lemma}

  \begin{Theorem}\cite{M1}\label{DNR}
  Let $M,N$ be two finite dimensional $H$-modules. Then, $M^*$ is an $H$-module with module structure defined as $(hf)(m)=f(S(h)m)$, for any $h\in H,\ m\in M,\ f\in M^*$, and
  \begin{eqnarray*}
&&M^*\otimes N^*\cong (N\otimes M)^*,\\
  &&M\otimes N^*\cong \Hom (N,M),\\
  &&\Hom_H(M,N)\cong \Hom (M,N)^H,
  \end{eqnarray*}
  where for any $H$-module $M$, $M^H=\{m\in M|hm=\epsilon(h)m,\forall h\in H\}$.
  \end{Theorem}

  \begin{Definition}\cite{A,DNR,M}
  Let $H$ be a Hopf algebra, $A$ be an associative algebra. $A$ is called an $H$-module algebra, if:

  (i) $A$ is an $H$-module.

  (ii) $h(ab)=\sum(h_1a)(h_2b),$ for any $h\in H,a,b\in A,$ where $\Delta(h)=\sum h_1\otimes h_2.$
  \end{Definition}

\section{Hopf actions on vertex operator algebras}

In this section, we first  recall  basic facts on vertex operator algebras from \cite{FLM, LL}. We then define so called $H$-module vertex operator algebra $V$ and 
discuss some consequences. We also present a Schur-Weyl type decomposition of $V$ as a module for the pair $(H,V^H)$, where $V^H$ is the $H$-invariants
of $V.$ These results generalize those given in \cite{DLM3}  in the case $H$ is the group algebra of a compact Lie group which acts on $V$ continuously.

  \begin{Definition}
  Let $V=\oplus_{n\in \mathbb{Z}}V_n$ be a $\mathbb{Z}$-graded vector space with $V_n=0$, $0\gg n$, $\dim V_n<\infty,\forall n\in \mathbb{Z}$, and $\textbf{1}\in V_0,\omega\in V_2,Y(\cdot,z):V\otimes V\rightarrow \End(V)[[z,z^{-1}]],u\otimes v\mapsto Y(u,z)v=\sum_{n\in\mathbb{Z}}u_nvz^{-n-1}$, where $u_n\in \End(V)$. Then, $V=(V,Y,\textbf{1},\omega)$ is called a vertex operator algebra(VOA for short) if the following hold:

  (i) For $ u,v\in V,$ $u_nv=0$, if $n\gg0$.

  (ii) $Y(\textbf{1},z)v=v,\lim_{z\rightarrow0}Y(v,z)\textbf{1}=v$, for $ v\in V$.

  (iii) Write $Y(\omega,z)=\sum_{n\in\mathbb{Z}}\omega_nz^{-n-1}=\sum_{n\in\mathbb{Z}}L(n)z^{-n-2}$, then
  \begin{eqnarray*}
  &&V_n=\{v\in V|L(0)v=nv\},~Y(L(-1)v,z)=\frac{d}{dz}Y(v,z),\\
  &&[L(n),L(m)]=(n-m)L(n+m)+\delta_{n+m,0}\frac{n^3-n}{12}c,
  \end{eqnarray*}
  where $c\in\mathbb{C}$ is called central charge of $V$. For $v\in V_n$, $v$ is said to be homogeneous and the weight $\wt v$ of $v$ is defined to be $n$.

  (iv) For $ u,v\in V$, we have
  \begin{eqnarray*}
  & &z_0^{-1}\delta(\frac{z_1-z_2}{z_0})Y(u,z_1)Y(v,z_2)-z_0^{-1}\delta(\frac{-z_2+z_1}{z_0})Y(v,z_2)Y(u,z_1)\\
  & &\ \ \ \ \ =z_1^{-1}\delta(\frac{z_2+z_0}{z_1})Y(Y(u,z_0)v,z_2),
  \end{eqnarray*}
  where $\delta(z)=\sum_{n\in\mathbb{Z}}z^n$, and $(z_1+z_2)^n=\sum_{i\geq0}\tbinom{n}{i}z_1^iz_2^{n-i}.$
  \end{Definition}

  Let $V$ be a vertex operator algebra. A {\em weak  $V$-module} $M$ is a vector space equipped
with a linear map
\begin{align*}
Y_{M}:V&\to (\End M)[[z, z^{-1}]],\\
v&\mapsto Y_{M}(v,z)=\sum_{n\in\Z}v_nz^{-n-1},\,v_n\in \End M,
\end{align*}
satisfying the following conditions: For any $u\in V,\ v\in V,\ w\in M$ and $n\in \Z$,
\begin{align*}
&\ \ \ \ \ \ \ \ \ \ \ \ \ \ \ \ \ \ \ \ \ \ \ \ \ \ u_nw=0 \text{ for } n\gg0;\\
&\ \ \ \ \ \ \ \ \ \ \ \ \ \ \ \ \ \ \ \ \ \ \ \ \ \ Y_M(\1, x)=\Id_M;\\
\begin{split}
&x_{0}^{-1}\delta\left(\frac{z_{1}-z_{2}}{z_{0}}\right)Y_{M}(u,z_{1})Y_M(v,z_{2})-z_{0}^{-1}\delta\left(
\frac{z_{2}-z_{1}}{-z_{0}}\right)Y_M(v,z_{2})Y_M(u,z_{1})\\
&\quad=x_{2}^{-1}\delta\left(\frac{z_{1}-z_{0}}{z_{2}}\right)Y_M(Y(u,z_{0})v,z_{2}).
\end{split}
\end{align*}

A weak
 $V$-module  $M$ is called an \textit{admissible $V$-module} if $M$ has a $\Z_{\geq
0}$-gradation $M=\bigoplus_{n\in\Z_{\geq 0}}M(n)$ such
that
\begin{align*}\label{AD1}
a_mM(n)\subset M(\wt{a}+n-m-1)
\end{align*}
for any homogeneous $a\in V$ and $m,\,n\in\Z$.
An admissible $V$-module $M$ is said to be
\textit{irreducible} if $M$ has no non-trivial admissible
$V$-submodule. When an admissible $V$-module $M$ is
direct sum of irreducible admissible submodules, $M$ is called
\textit{completely reducible}.

A vertex operator algebra $V$ is said to be \textit{rational} if
any  admissible $V$-module is completely reducible.

A {\em  $V$-module} is a weak $V$-module $M$ which carries a $\C$-grading induced by the spectrum of $L(0)$, that is  $M=\bigoplus_{\lambda\in\C}
M_{\lambda}$ where
$M_\lambda=\{w\in M|L(0)w=\lambda w\}$. Moreover one requires that $M_\lambda$ is
finite dimensional and for fixed $\lambda\in\C$, $M_{\lambda+n}=0$
for sufficiently small integer $n$.

The following results will be useful in this paper.

  \begin{Lemma}\cite{DM,L}\label{span}
  Let $V$ be a simple VOA. Then $V$ is spanned by $u_nv,u\in V, n\in \mathbb{Z}$ for any fixed nonzero $v\in V.$
  \end{Lemma}

  \begin{Lemma}\cite{DM}\label{2}
  Let $V$ be a simple VOA, and $(M,Y_M)$ a simple $V$-module. Let $u^1,u^2,\cdots,u^s\in V,m^1,m^2,\cdots,m^s\in M$ such that $u^1,u^2,\cdots,u^s\in V$ are not zero and $m^1,m^2,\cdots,m^s$ are linearly independent. Then, $$\sum_{i=1}^sY_M(u^i,z)m^i\neq0.$$
  \end{Lemma}

  \begin{Remark}
  For any VOA, the skew-symmetry identity $$Y(u,z)v=e^{zL(-1)}Y(v,-z)u$$ is always true \cite{LL}.
  \end{Remark}

  Now, we define the main objects of this paper.

  \begin{Definition}\label{D}
  Let $(V,Y,\textbf{1},\omega)$ be a vertex operator algebra, and $H$ be a Hopf algebra, we say $V$ is an $H$-module VOA if

  (i) $V$ is an $H$-module.

  (ii) $h\omega=\epsilon(h)\omega$, for $ h\in H.$

  (iii) $hY(u,z)v=\sum Y(h_1u,z)h_2v,$ for $ h\in H,~u,v\in V,$
where $\Delta(h)=\sum h_1\otimes h_2$ is Sweedler's notation.
  \end{Definition}

Both $\Aut(V)$ and $\Der(V)$ can be understood in terms of Hopf algebra actions. Let $G=\Aut(V).$ Then the group algebra $\C[G]$ is a Hopf algebra such that $\Delta(g)=g\otimes g$ and $\epsilon(g)=1$ for $g\in G.$ Condition (iii) in this case becomes $gY(u,z)v=Y(gu,z)gv.$ This together with (ii) gives the definition of an automorphism. Similarly, $U(\Der(V))$ is a Hopf algebra and $V$ is a $U(\Der(V))$-module vertex operator algebra. So the notion of $H$-module vertex operator algebra unifies and generalizes both notions of automorphism and derivation of a vertex operator algebra. 

Here are some consequences of the definition:
\begin{Lemma}\label{V^H} Let $V$ be a simple  $H$-module VOA. Then

(i) For  any $h\in H,$  $hY(\omega,z)=Y(\omega,z)h.$ In particular, $h$ preserves each homogeneous subspace $V_n$ and $V_n$ is an $H$-module,

(ii)  $V^H=\{v\in V|hv=\epsilon(h)v\}$ is a vertex operator subalgebra of $V.$ $V^H$ is called the $H$-invariants of $V.$
\end{Lemma} 

\proof (i) By the definition, for $ h\in H, v\in V$, we have
  \begin{eqnarray*}
  hY(\omega,z)v=\sum Y(h_1\omega,z)h_2v=\sum Y(\epsilon(h_1)\omega,z)h_2v=Y(\omega,z)hv.
\end{eqnarray*}
The rest of (i) is clear.

(ii) Obviously， $\omega\in V^H$, so $V^H\neq0$. From
  \begin{eqnarray*}
  Y(\epsilon(h)\textbf{1},z)\omega&=&\epsilon(h)\omega=h\omega=hY(\textbf{1},z)\omega=\sum Y(h_1\textbf{1},z)h_2\omega\\
  &=&\sum Y(h_1\epsilon(h_2)\textbf{1},z)\omega=Y(h\textbf{1},z)\omega,
  \end{eqnarray*}
  we have $Y(\epsilon(h)\textbf{1}-h\textbf{1},z)\omega=0.$ Since $V$ is simple and $\omega\ne 0$ , $h\textbf{1}=\epsilon(h)\textbf{1}$  by Lemma \ref{2}. That is, $\textbf{1}\in V^H.$

It remains to show $v_nu\in V^H$ for $ n\in \mathbb{Z}, u,v\in V^H.$ But this follows from
  \begin{eqnarray*}
  hY(v,z)u=\sum Y(h_1v,z)h_2u=\sum \epsilon(h_1)\epsilon(h_2)Y(v,z)u=\epsilon(h)Y(v,z)u
  \end{eqnarray*}
immediately for $h\in H.$  $\hfill\Box$

 We now assume that $H$-module VOA $V$ is a simple VOA and  a semisimple $H$-module. Let $\{M_i|i\in I\}$ be the set of inequivalent simple $H$-modules occurring in $V.$ Since $H$ preserves each homogeneous subspace $V_n$ of $V$ by Lemma \ref{V^H} (i), we see that
each $M_i$ is a finite dimensional $H$-module. Let $0\in I$ such that $M_0=\C$ is the simple $H$-module such that $h$ acts 
as $\epsilon(h)$ for $h\in H.$ Let $V^i$ be the sum of simple $H$-submodules of $V$ isomorphic to $M_i.$ Then $V^0=V^H.$ 
 Then as an $H$-module, $V$ has the following decomposition:
$$V=\oplus_{i\in I}V^i=\oplus_{i\in I} M_i\otimes V_{M_i},$$
 where $V_{M_i}$ is the multiplicity space of $M_i$ appearing in $V.$ Note that $V_{M_0}=V^H.$  

In fact, $V_{M_i}$ can be realized as a subsapce of $V^i.$ Note that $V^{i}_n=V_n\cap V^{i}$ is the direct sum of copies of $M^i.$ Fix any nonzero vector $w$ of $M_i.$  We can set
$$V_{M_i}=\{f(w)|f\in \Hom_H(M_i,V)\}.$$ 
For example, we can take $w$ be a highest weight vector of $M_i$ in the following sense. Since the action of $H$ on $M_i$ induces an algebra epimorphism $H\rightarrow \End(M_i)$, we define highest-weight vectors in $M_i$ as follows: Fix a basis of $M_i$ with $E^i_{pq}$ be corresponding standard basis of $\End(M_i)$. Then $u\in M_i$ is a highest-weight vector if $E^i_{pq}u=0$ whenever $q>p$ or $p=q,p>1$; and if $E^i_{11}u=u$. Take $w=u,$ then $V_{M_i}$ is the space of  the highest weight vectors in $V^{i}$. 

\begin{Lemma}\label{V^H1} The actions of $H$ and $V^H$ commute. In particuar, each $V_{M_i}$ is a $V^H$-module. In other words, $H$ acts on the first tensor factor of $V^{i}=M_i\otimes V_{M_i}$, and $V^H$ acts on the second tensor factor.

\end{Lemma}

\proof Let $u\in V^H, v\in V$ and $h\in H.$ Then 
\begin{eqnarray*}
hY(u,z)v &=& \sum Y(h_1u,z)h_2v\\
&=&\sum Y(\epsilon(h_1)u,z)h_2v\\
&=&\sum Y(u,z)\epsilon(h_1)h_2v\\
& =&Y(u,z)hv
\end{eqnarray*}
 and $Y(u,z)$ and $h$ commute. As a result, $V^H$ acts on the multiplicity space $V_{M_i}.$
 $\hfill\Box$

 \begin{Theorem}\label{T} Let $V,$ $H,$ $V_{M_i}$ be as before. Then we have the following Schur-Weyl type duality:

  (i) $V^H$ is a simple VOA and each  $V_{M_i}$ is a simple $V^H$-module.

  (ii) $V_{M_i}$ and $V_{M_j}$ are isomorphic $V^H$-modules if and only if $i= j.$
\end{Theorem}

  \proof The main tool we use to prove the theorem is the $A_n(V)$-theory developed in \cite{DLM1}.  We first recall from \cite{DLM1} the definition of
$A_n(V)$ with $n\geq 0.$  For homogeneous $u,v\in V,$ define two bilinear operations $\ast_n,\circ_n$ on $V$  as follows:
  \begin{eqnarray*}
  &&u\ast_n v=\sum^n_{m=0}(-1)^m\tbinom{m+n}{m}\Res_zY(u,z)v\frac{(1+z)^{\wt u+n}}{z^{n+m+1}},\\
  &&u\circ_n v=\Res_zY(u,z)v\frac{(1+z)^{\wt u+n}}{z^{2n+2}}.
  \end{eqnarray*}
Set $O_n(V)=\Span\{u\circ_n v,~L(-1)u+L(0)u|\forall u,v\in V\},$ and $A_n(V)=V/O_n(V)$.  Here we present some results from \cite{DLM1} which will be used in the proof: (1) $A_n(V)$ is an associative algebra, (2) The identity map on $V$ induces an onto algebra homomorphism from $A_n(V)$ to $A_m(V)$ for $m\leq n,$ (3) If $M=\oplus_{m\geq 0}M(m)$ is an admissible $V$-module then each $M(s)$ is an $A_n(V)$-module
for $s\leq n$ where $v+O_n(V)$ acts as $o(v)$ and $o(v)=v_{\wt v-1}$ for homogeneous $v.$ Moreover, $M$ is irreducible if and only if $M(n)$ is a simple $A_n(V)$-module for all $n\geq 0.$

We now divide the proof into several parts:

(1)  $A_n(V)$ is an $H$-module algebra for each $n\geq0.$ Note that  $H$ preserves the $L(0)$-grading of $V.$  Let  $u,v$ be homogeneous. Then
  \begin{eqnarray*}
  &&h(u\circ_n v)=\Res_zhY(u,z)v\frac{(1+z)^{\wt~u+n}}{z^{2n+2}}\\
  &=&\sum \Res_zY(h_1u,z)h_2v\frac{(1+z)^{\wt~u+n}}{z^{2n+2}}=\sum h_1u\circ_n h_2v.
  \end{eqnarray*}
  Clearly, $h(L(-1)u+L(0)u)=(L(-1)u+L(0)(hu).$ So, $H$ preserves $O_n(V)$ and $A_n(V)$ is an $H$-module. Similarly, $h(u\ast_n v)=\sum h_1u\ast_n h_2v.$ So $A_n(V)$ is an $H$-module algebra. Moreover, as a quotient module of a semisimple $H$-module, $A_n(V)$ is a semisimple $H$-module.

(2) The $\End (V)$ is an $H$-module via conjugation, i.e. $$(h\cdot f)u=\sum h_1f(S(h_2)u)$$
for $h\in H,$ $f\in \End(V)$ and $u\in V.$ In fact, each $\End(V_n)$ is a finite dimensional $H$-module.

(3) Fix $n\geq 0,$ the map $\sigma_n: u\mapsto o(u)|_{V_n}$ is an onto $H$-module homomorphism. In particular,
$\End(V_n)$ is a semisimple $H$-module. Clearly,  the linear map  is a well-defined. For  $v\in V, u\in V_n, h\in H$, we have
  \begin{eqnarray*}
 (h\cdot o(v))u&=&\sum h_1o(v)S(h_2)u\\
&=&\sum o(h_1v)(h_2S(h_3)u)\\
&=&\sum o(h_1v)(\epsilon(h_2)u)\\
  &=&\sum o(h_1\epsilon(h_2)v)u\\
&=&o(hv)u.
  \end{eqnarray*}
  Thus $\sigma_n$  is an $H$-module homomorphism from $A_n(V)$ to $\End(V_n).$ As $V_n$ is a simple $A_n(V)$-module, the map is onto. It is clear now that $\sigma(A_n(V)^H)=\End(V_n)^H.$ That is, $\{o(u)|_{V_n}|u\in V^H\}=\{f\in \End(V_n)|h\cdot f=\epsilon(h)f, \ h\in H\}.$

(4) $\End(V_n)^H=\End_H(V_n).$  This is true in general (See \cite{M1} or Theorem \ref{DNR}). In the current situation, we present a very short proof. If $f\in \End_H(V_n)$, then
$$h\cdot f=\sum h_1fS(h_2)=\sum fh_1S(h_2)=\epsilon(h)f.$$
 That is, $f\in \End(V_n)^H.$ 
Now let $f\in \End(V_n)^H,$ there exists $u\in V^H$ such that $f=o(u)|_{V_n}.$ By Lemma \ref{V^H1}, $ hf=ho(u)=o(u)h$ on $V_n$ for all $h\in H.$ Thus $f\in \End_H(V_n).$

  (5) $\End(V_n)$ is an $H$-module algebra.  From (3), $\sigma_n$ is onto.  is surjective. Let $u, v\in V, $ $w\in V_{n},h\in H$. Then 
  \begin{eqnarray*}
  (h\cdot o(u)o(v))w&=&\sum h_1o(u*_nv)(S(h_2)w))\\
  &=&\sum o(h_1(u*_nv))(h_2S(h_3)w)\\
&=&\sum o((h_1u)*_n(h_2v))\epsilon(h_3)w\\
  &=&\sum o(h_1u)o(\epsilon(h_3)h_2v)(w)\\
&=&\sum o(h_1u)o(h_2v)w.
  \end{eqnarray*}
  So, we have $h\cdot(o(u)o(v))=\sum o(h_1u)o(h_2v)$ on $V_n.$ From (3), we know that $h\cdot o(u)=o(hu).$ As a result, $\sum o(h_1u)o(h_2v)=
(h_1\cdot o(u))(h_2\cdot o(v)).$
 
The rest of the proof is similar to that given in \cite{DLM3}.

(6)  For $i\in I,$ $V_{M_i}$ is a simple $V^H$-module. For convenience we can asume that
$V=\oplus_{n\geq 0}V_n.$ Note that $V_{M_i}=\oplus_{n\geq 0}(V_{M_i})_n$, where
$(V_{M_i})_n=V_n\cap V_{M_i}.$ From the $A_n(V)$ theory, it is good enough to prove that
each $(V_{M_i})_n\ne 0$ is a simple $A_n(V^H)$-module.

 It is sufficient to show that if  $(V_{M_i})_n\ne 0$, then  $(V_{M_i})_n$ is generated
 by any nonzero vector of $(V_{M_i})_n$ as $A_n(V^H)$-module. Let $u,v\in (V_{M_i})_n$ to be  linearly independent. Then as an $H$-module we can write $V_n$ as direct sum of three submodules $M_i\otimes u\oplus M_i\otimes v\oplus W$. Define a map $\alpha\in \End(V_n)$ such that $$\alpha(x\otimes u+y\otimes v+w)=y\otimes u+x\otimes v+w,$$ for $x,y\in M_i, w\in W$. Clearly, $\alpha\in \End_H(V_n).$ Hence, $\alpha\in\End(V_n)^H$ by (4). It follows from (3) there exists $a\in V^H$ such that
$\alpha=o(a)$ on $V_n.$ Then $o(a)u=v$ and $(V_{M_i})_n$ can be generated by any nonzero vector $u$ as $A_n(V^H)$-module.

(7) If $i,j$ are different, then $V_{M_i}$ and $V_{M_j}$ are not isomorphic $V^H$-modules.
 Pick up  $n$ such that $V^i_n\neq0.$  Then $V_n$ is direct sum of three $H$-modules $V_n=V_n^i\oplus V_n^j\oplus W$ for a suitable $H$-submodule $W$ of $V_n$. Define $\beta\in \End(V_n)$ such that it is identity on $V_{n}^i$ and zero on $V_{n}^j$ and $W$. Clearly $\beta\in \End_H(V_n)=\End(V_n)^H$. Again by (3) there exists $b\in V^H$ such that 
$\beta=o(b).$  Thus there is no $V^H$-module homomorphism between $V_{M_i}$ and $V_{M_j}$,
and  $V_{M_i}$ and $V_{M_j}$ are isomorphic $V^H$-modules if and only if $i= j.$ $\hfill\Box$

\begin{Remark}\label{Cjt} In the case that $H$ is a group algebra $\C[G]$ where $G$ is a compact Lie group which acts on $V$ continuously, Theorem \ref{T} was obtained previously in \cite{DLM3}. Moreover, if $G$ is a finite group and $V^G$ is rational, $C_2$-cofinite,
then $\qdim_{V^G}V_{M_i}=\dim M_i$ \cite{DJX}, where $\qdim_{V^G}V_{M_i}$ is the quantum
dimension of $V_{M_i}$ over $V^G$ \cite{DLM3}.  Motivated by this result, we conjecture
that if $V$ is a simple VOA and a semisimple $H$-module VOA then $\qdim_{V^H}V_{M_i}=\dim M_i$ for all $i\in I.$ But even in the case $H=\C[G]$ and $G$ is a compact Lie group in general, we do not know how to establish this conjecture.
\end{Remark}

\section{Intertwining operators}

We are working in the settings of Section 3. In particular, $V$ is a simple and semisimple $H$-module VOA where
$H$ is a Hopf algebra. In this section, we discuss how the fusion rules $m_{ij}^l$ between $H$-modules $M_i,M_j$ and $M_l$  has empact on the fusion rules  $N_{V_{M_i}\,V_{M_j}}^{V_{M_l}}$ between nonzero $V^H$-modules $V_{M_i},V_{M_j}$, and $V_{M_l}.$

We first define the fusion rules $m_{ij}^l$ and $N_{M_i\,M_j}^{M_l}$ for $i,j,l\in I.$ Let $H$ be a Hopf algebra, and $M_i,M_j,M_l$ be three finite dimensional $H$-modules, we define $m_{M_i\,M_j}^{M_l}=\dim \Hom_H(M_l,M_i\otimes M_j)$.   For short we set $m_{ij}^l=m_{M_i\,M_j}^{M_l}.$ Clearly,  $m_{ij}^l$ are finite as 
$M_i,M_j, M_l$ are finite dimensional.

 The  fusion rules \cite{FHL} among $V$-modules are more complicated. 
  
\begin{Definition}
  Let $U^1,U^2,U^3$ be three weak $V$-modules. An intertwining operator $\mathcal {Y}(\cdot,z)$ of type $\tbinom{U^3}{U^1~~U^2}$ is a linear map
   \begin{eqnarray*}
  \mathcal {Y}(\cdot,z):~U^1&\rightarrow&~\Hom(U^2,U^3)\{z\},\\
  v^1&\mapsto& \mathcal {Y}(v^1,z)=\sum_{n\in\mathbb{C}}v_n^1z^{-n-1},
  \end{eqnarray*}
  where $v_n^1\in \Hom(U^2,U^3)$, satisfying the following conditions:

  (i) For any $v^1\in U^1,v^2\in U^2,\lambda\in\mathbb{C}$, $v^1_{n+\lambda}v^2=0$ for $n\in\mathbb{Z}$ sufficiently large.

  (ii) For any $v^1\in U^1$, $\frac{d}{dz}\mathcal {Y}(v^1,z)=\mathcal {Y}(L(-1)v^1,z)$.

  (iii) For any $a\in V, v^1\in U^1$,
  \begin{eqnarray*}
  z_0^{-1}\delta(\frac{z_1-z_2}{z_0})Y_{U^3}(a,z_1)\mathcal {Y}(v^1,z_2)-z_0^{-1}\delta(\frac{-z_2+z_1}{z_0})\mathcal {Y}(v^1,z_2)Y_{U^2}(a,z_1)\\
  =z_1^{-1}\delta(\frac{z_2+z_0}{z_1})\mathcal {Y}(Y_{U^1}(a,z_0)v^1,z_2),
  \end{eqnarray*}
  \end{Definition}

  Obviously, all the intertwining operators of type $\tbinom{U^3}{U^1~~U^2}$ form a vector space, denoted by $I_V\tbinom{U^3}{U^1~~U^2}.$  The number $N_{U^1,U^2}^{U^3}=\dim I_V\tbinom{U^3}{U^1~~U^2}$.  $N_{U^1\,U^2}^{U^3}$ is called the fusion rules of type $\tbinom{U^3}{U^1~~U^2}$.

Then, we have the following proposition which is similar to Theorem 2 in \cite{T}.

  \begin{Proposition}\label{0}
  Let $V$ be a simple $H$-module VOA, and semisimple as $H$-module with decomposition $V=\oplus _{i\in I}M_i\otimes V_{M_i}$ as in Section 3.  Then, $N_{V_{M_i}\,V_{M_j}}^{V_{M_l}}\geq m_{ij}^l$.
  \end{Proposition}
  \proof Since $N_{V_{M_i}\,V_{M_j}}^{V_{M_l}}=\dim I_{V^H}\tbinom{V_{M_l}}{V_{M_i}~~V_{M_j}}$, and $m_{ij}^l=\dim \Hom_H(M_l,M_i\otimes M_j).$ We just need to find an injective linear map from $\Hom_H(M_l,M_i\otimes M_j)$ to $I_{V^H}\tbinom{V_{M_l}}{V_{M_i}~~V_{M_j}}$.

  Since $I_{V}\tbinom{V}{V~~V}=\mathbb{C}Y$, we consider $\mathbb{C}Y, V_{M_i},V_{M_j}$ as trivial $H$-modules. For any $ v^{10}\in V_{M_i},v^{20}\in V_{M_j}$, set $$\mathcal {I}_{12}(v^{10},v^{20})=\Span\{Y\otimes(w^1\otimes v^{10})\otimes(w^2\otimes v^{20})|w^1\in M_i,w^2\in M_j\}.$$
  Then, $\mathcal {I}_{12}(v^{10},v^{20})$ is an $H$-module and isomorphic to $M_i\otimes M_j$.

  Let $F\in \Hom_H(M_l,\mathcal {I}_{12}(v^{10},v^{20}))$, we shall define $\Phi(F)\in I_{V^H}\tbinom{V_{M_l}}{V_{M_i}~~V_{M_j}}$ as follows.

  For $w^{31}\in M_l$, we set $F(w^{31})=\sum_sY\otimes(w^{1s}\otimes v^{10})\otimes(w^{2s}\otimes v^{20})$ for some $w^{1s}\in M_i,w^{2s}\in M_j$. Now, for $v^!\in V_{M_i},v^2\in V_{M_j}$, we define $$H(v^1,v^2)(w^{31})\triangleq\sum_sY\otimes(w^{1s}\otimes v^{1})\otimes(w^{2s}\otimes v^{2}).$$
  It is clear that $H(v^1,v^2)\in \Hom_H(M_l,\mathcal {I}_{12}(v^{1},v^{2}))$.

  We define a linear map $\Psi$ as $$\Psi(\sum_sY\otimes(w^{1s}\otimes v^{1})\otimes(w^{2s}\otimes v^{2}))\triangleq\sum_sY({w^{1s}\otimes v^{1},z})(w^{2s}\otimes v^{2}).$$
  Then, the image of $\Psi$ is indeed in $V^l\{z\}$, where $V^l=M_l\otimes V_{M_l}$. Consider $V^l\{z\}$ as $H$-module naturally induced from $H$-module $V^l$. Then, $\Psi$ is an $H$-module homomorphism. So $\Psi(H(v^1,v^2)(M_l))$ is an $H$-submodule of $V^l\{z\}$

  Let $w^{31},\cdots,w^{3\dim M^l}$ be a basis of $M_l$, there exists $h\in H$ such that $hw^{3s}=\delta_{1s}w^{31}$. Write $\Psi(H(v^1,v^2)(w^{31}))=\sum_sw^{3s}\otimes p^s$, where $p^s\in V_{M_l}\{z\}$. Then
  \begin{eqnarray*}
  \Psi(H(v^1,v^2)(w^{31}))&=&\Psi(H(v^1,v^2)(hw^{31}))=h\Psi(H(v^1,v^2)(w^{31}))\\
  &=&h\sum_sw^{3s}\otimes p^s=\sum_shw^{3s}\otimes p^s=w^{31}\otimes p^1.
  \end{eqnarray*}
  We hence have an unique $\Phi(F)(v^1,z)v^2\in V_{M_l}\{z\}$ such that$$w^{31}\otimes\Phi(F)(v^1,z)v^2=\Psi(H(v^1,v^2)(w^{31})).$$
  Since $Y$ is an intertwining operator, we have $\Phi(F)\in I_{V^H}\tbinom{V_{M_l}}{V_{M_i}~~V_{M_j}}.$

  Now, we show that $\Phi$ is injective. Suppose $\Phi(F)=0,$ then $$w^{31}\otimes\Phi(F)(v^1,z)v^2=\Psi(H(v^{10},v^{20})(w^{31}))=\Psi(F(w^{31}))=0.$$
  By Lemma \ref{2}, we have $F(w^{31})=0$. Since $F$ is an $H$-module homomorphism and $M_l$ is irreducible, we have $F=0.$$\hfill\Box$

\begin{Remark} In the case $H$ is a group algebra of a finite group, above Proposition was obtained previously in \cite{T}. 
In fact, the inequality becomes equality in \cite{DRX} in this case with the help of the quantum dimension identity \cite{DJX}. We certainly believe that $N_{M_i,M_j}^{M_k}=m_{ij}^k$ for any Hopf algebra $H.$ But we do not have a clue on how to prove it in this paper.
\end{Remark}

  \section{Faithful Hopf actions}

 In this section, we always assume that $H$ is a finite dimensional  Hopf algebra and $V$ is an $H$-module simple VOA.

\begin{Lemma}\label{l5.1} Suppose $V$ is an $H$-module VOA.  Then for any two finite dimensional $H$-submodule $M,N$ of $V,$  $\varphi:M\otimes N\rightarrow N\otimes M,u\otimes v\mapsto v\otimes u$ is an $H$-module isomorphism. In particular, \begin{eqnarray}\label{Cocm}
\sum h_1v\otimes h_2u=\sum h_2v\otimes h_1u
\end{eqnarray}
 for any $h\in H,u\in M,v\in N$.
  \end{Lemma}

\proof Let $X$ be the linear subspace of $V[[z,z^{-1}]]$ spanned by $Y(u,z)v$  and $Z$ be the linear subspace of $V[[z,z^{-1}]]$ spanned by $Y(v,z)u$ for $u\in M$ and $v\in N$. Regard $V[[z,z^{-1}]]$ as an $H$-module naturally induced from the  $H$-module  structure on $V$. By the definition of $H$-module VOA, it is easy to see that both $X$ and $Z$ are $H$-submodules of $V[[z,z^{-1}]]$. Define $Y':M\otimes N\rightarrow X$ as $Y'(u\otimes v)=Y(u,z)v$. Since $h(Y(u,z)v)=\sum Y(h_1u,z)h_2v$, $Y'$ is an $H$-module homomorphism. By Lemma \ref{2}, $Y'$ is an $H$-module isomorphism. Similarly, $N\otimes M$ and $Z$ are isomorphic $H$-modules.

  Define $f:X\rightarrow Z,Y(u,z)v\mapsto Y(v,z)u$, again by Lemma \ref{2}, $f$ is a linear isomorphism. Since $Y(u,z)v=e^{zL(-1)}Y(v,-z)u$ for $u\in M,v\in N$, we have
  \begin{eqnarray}
  hf(Y(u,z)v)&=&h(Y(v,z)u)\nonumber\\
  &=&h(e^{zL(-1)}Y(u,-z)v)\nonumber\\
  &=&\sum e^{zL(-1)}Y(h_1u,-z)h_2v\nonumber\\
  &=&\sum Y(h_2v,z)h_1u\nonumber\\
  &=&f(\sum Y(h_1u,z)h_2v)\nonumber\\
  &=&f(hY(u,z)v).\label{.}
  \end{eqnarray}
 Here we use that fact that  $h$ commutes with $Y(\omega,z)$ from  Lemma \ref{V^H} (i). So, $f$ is an $H$-module isomorphism. As a result, $\varphi$ is an $H$-module isomorphism. The isomorphism implies (\ref{Cocm}).
$\hfill\Box$

\begin{Theorem}\label{t5.2}
In Lemma \ref{l5.1}, A faithful action of a finite dimensional Hopf algebra $H$ on a simple VOA is a group action.
\end{Theorem}
\proof By Lemma \ref{DNR}, it is enough  to show $H$ is cocommutative. Since $V$ is a faithful $H$-module there 
 exists a positive integer $n$ such that $V_{\leq n}=\oplus_{m\leq n}V_m$ is a faithful $H$-module.
Consider $V_{\leq n}\otimes V_{\leq n}$ as $H\otimes H$-module with module structure $$(h\otimes g)\cdot(u\otimes v)=hu\otimes gv$$ for any $h,g\in H,u,v\in V_{\leq n} $. Then, $V_{\leq n}\otimes V_{\leq n}$ is a faithful $H\otimes H$-module. Consider $\sum h_1\otimes h_2$ and $\sum h_2\otimes h_1$ as elements of $H\otimes H$. From identity (\ref{Cocm}) with $M=N=V_{\leq n}$  we have
$$
\sum h_1v\otimes h_2u=\sum h_2v\otimes h_1u
$$
 for any $h\in H, u, v\in V_{\leq n}.$  we have  $\sum h_1\otimes h_2=\sum h_2\otimes h_1$. 
Hence, $H$ is cocommutative, and thus a group algebra.$\hfill\Box$

\section{Inner faithful Hopf actions}

 In this section, we always assume that $H$ is a finite dimensional Hopf algebra. Also assume the $H$-module VOA is a simple VOA and semisimple $H$-module.  Let $\{M_i|i\in I\}$ be the set of inequivalent simple $H$-modules occurring in $V$ with $0\in I$ and 
 $M_0$  the trivial module. Let $\psi_i:H\rightarrow \End(M_i)$ be the $i^{th}$ simple representation. The character $\chi_i\in H^*$ is defined by $\chi_i(h)=Tr_{M_i}\psi_i(h)$, then $\chi_0=\epsilon$. 
Note that  $M_i^*$ is also an irreducible left $H$-module such that  $f\in M_i^*,\ m\in M_i$, $$(hf)(m)=f(S(h)m)$$ for any $h\in H$. We sometimes write $M_{Si}$ for $M_{i}^*$ and $\chi_{i^*}$ for $\chi_{M_i^*}=\chi_{Si}$. 

   From discussion in Section 3, we can write the decomposition of $V$ as $$V=\oplus_{i\in I}V^i=\oplus_{i\in I}M_i\otimes V_{\chi_i}.$$

In the  following two definitions, $H$ can be any Hopf algebras, not necessary finite dimensional Hopf algebras.

  \begin{Definition}\cite{A,DNR,M}
  Let $H$ be a Hopf algebra and $I$ be a subspace of $H$. If

  (i) $I$ is an ideal of $H$ as associative algebra, i.e. $IH\subseteq I,HI\subseteq I.$

  (ii) $I$ is a coideal of $H$ as coassociative coalgebra, i.e. $\Delta(I)\subseteq H\otimes I+I\otimes H$, and $\epsilon(I)=0$.

  (iii) $I$ is stable under the antipode $S$ of $H$, i.e. $S(I)\subseteq I.$

  Then, $I$ is called a Hopf ideal. If $I$ is a Hopf ideal, then $H/I$ is the quotient Hopf algebra.
  \end{Definition}
  \begin{Definition}\cite{EW1}
  Given a left $H$-module $M$, we say that $M$ is inner faithful if $IM\neq0$ for any nonzero Hopf ideal $I.$ The action of $H$ on a VOA $V$ is inner faithful if the left $H$-module $V$ is inner faithful.
  \end{Definition}
  
 \begin{Lemma}\label{h}
  Let $H$ be a  Hopf algebra, $V$ be an $H$-module simple VOA , and $K$ be the kernel of the action of $H$ on $V$. Then, $\epsilon (K)=0$, and $S(K)\subseteq K$.
  \end{Lemma}
 \proof First, we prove $\epsilon(K)=0$. For any $ k\in K$, we have $$0=k\omega=\epsilon(k)\omega,$$
  so $\epsilon(K)=0$.

  Now, we prove $K$ is stable under $S$.  From Theorem \ref{DNR}, Lemma \ref{l5.1} we have $$(M_{{S_j}}\otimes M_i)^H=(M_i\otimes M_{{S_j}})^H=(M_i\otimes M_j^*)^H=\Hom (M_j,M_i)^H=\Hom _H(M_j,M_i)=\delta_{ij}\mathbb{C}.$$
If $V_{\chi_i}\ne 0$ then for any fixed nonzero $v\in V^i$, $V$ is spanned by $\{u_nv|u\in V,n\in\mathbb{Z}\}$
by Lemma \ref{span} . Since $\omega\in V^H$ we see that $V_{\chi_{Si}}\neq0.$ So, there is $j$, such that $M_j=M_{Si}$, and $V_{\chi_{j}}\neq0.$ Thus, $S$ induces a permutation on the index set $I$.

To show that $S(k)\subset K$, it is good enough to  establish that $S(K)M_i=0$ for $i\in I.$ Now fix $i\in I.$ Consider $M_{S^ri}$ for $r\geq 0.$ From the discussion above, we know that $M_{S^ri}$ occurs in $V$ for all $r$ and $M_{S^ri}=M_{S^{r-1}i}^*.$  Since
$H$ is finite dimensional, there is $n\geq 1$ such that $S^n={\rm Id}$ by Lemma \ref{A,DNR}.  So $M_{S^ni}=M_i.$ 
Let $k\in K$ and $f\in M_i.$ Since $M_i=M_{S^{n-1}i}^*$ we  need to show that $S(k)f(M_{S^{n-1}i})=f(S^2(k)M_{S^{n-1}i})=0.$
It suffices to show that $S^2(k)M_{S^{n-1}i}=0.$. Continuing in this way reduces to prove $S^n(k)M_{Si}=kM_{Si}=0.$  But this is clear. 
 $\hfill\Box$


\begin{Lemma}\label{l6.4}
  Let $A$ be a finite dimensional semisimple associative algebra and $M$ be a faithful $A$-module. Write $$A=\oplus_{i=1}^sM(n_i,\mathbb{C}),$$and let $E_{pq}^i,p,q=1,\cdots,n_i$ be the standard basis of each $M(n_i,\mathbb{C}).$ Fix $E_{p_0q_0}^{i_0}.$ Then  there is $m\in M$ such that $\{E_{pq_0}^{i_0}m|p=1,...,n_{i_0}\}$
is a linearly independent set and  all other $E_{pq}^im=0.$
  \end{Lemma}
  \proof Since $A$ is semisimple, $M$ is completely reducible. Without loss we can assume that $M=\oplus_{i=1}^s N_i$, where  $N_i$ is the simple  $M(n_i,\mathbb{C})$-module.

  Without loss we can assume that $q_0=1.$ Let $m=(1,0,\cdots,0)^T\in N_{i_0}.$ Then all nonzero $E^i_{pq}m$ are exactly $E^{i_0}_{p1}m$ for $p=1,...,n_{i_0}$. Clearly, these vectors are linearly independent. $\hfill\Box$

  \begin{Theorem}\label{wh}
 Let $H$ be a finite  dimensional Hopf algebra and $V$ be a simple VOA and  an inner faithful, semisimple  $H$-module VOA . Then, $H$ acts on $V$ faithfully. In particular, $H$ is a group algebra.
  \end{Theorem}
  \proof Let $K$ be the kernel of the given inner faithful action. 
Clearly  $K$ is an ideal of $H.$ Assume that $K$ is  not a Hopf ideal. Then $K$ is nonzero. By Lemma \ref{h}, we have $$\Delta (K)\nsubseteq H\otimes K+K\otimes H.$$

Denote by $J(H)$ the Jacobson radical of $H.$ Since $V$ is semisimple $H$-module, we have $J(H)\subset K$, and $H/K$ is a semisimple associative algebra acting on $V$ faithfully. Then there exists $n\geq 0$ such that  $V_{\leq n}=\oplus_{m\leq n}V_m$ is a faithful $H/K$-module.

 We have a decomposition of semisimple algebra $$H/K=\oplus_{i=1}^sM(n_i,\mathbb{C}),$$
  where $M(n_i,\mathbb{C})$ is the fully matrix algebra of degree $n_i$. Let $$\{E^i_{pq}|p,q=1,\cdots,n_i\}$$ be the standard basis of $M(n_i,\mathbb{C}),i=1,2,\cdots,s$. Let $e_{pq}^i$ be the elements in $H$ such that $\overline{e_{pq}^i}=e_{pq}^i+K=E_{pq}^i$, then we have $e_{pq}^iv=E_{pq}^iv$ for any $v\in V$. Let $k_1,\cdots,k_t$ be a basis of $K$, then $\{e_{pq}^i,k_j|i=1,2,\cdots,s,j=1,2,\cdots,t,p,q=1,\cdots,n_i\}$ is a basis of $H$.

 For $k\in K$, write $$\Delta(k)=\sum_{i=1}^s\sum_{p,q}h^i_{pq}\otimes e_{pq}^i+\sum_{j=1}^th_j\otimes k_j,$$ for some $k_j\in K,  h^i_{pq}, h_j\in H.$ Let $L$ be the linear subspace of $H$ spanned by $e_{pq}^i,\ i=1,2\cdots,s,\ p,q=1,\cdots, n_i$, then $H=K\oplus L$ as linear space. We can write $\Delta(k)$ as follows$$\Delta(k)=\sum_{j=1}^th_j\otimes k_j+\sum_{i=1}^s\sum_{p,q}(a^i_{pq}+b^i_{pq})\otimes e_{pq}^i,$$
  where $a^i_{pq}\in K,b^i_{pq}\in L$.

  Since  $\Delta (K)\nsubseteq H\otimes K+K\otimes H$, there is $k\in K$ such that  $b^{i_0}_{p_0q_0}\neq0$ for some $i_0<s$ and $1\leq p_0,q_0\leq n_{i_0}.$ Since $H/K$ is a semisimple associative algebra and $V_{\leq n}$ is a faithful $H/K$-module, it follows from Lemma  \ref{l6.4} that there is $v\in V$ such that $E_{pq_0}^{i_0}v$ for $p=1,...,n_{i_0}$ are linearly independent and other $E_{pq}^iv=0.$ Then for any $u\in V$, 
  \begin{eqnarray*}
  0&=&kY(u,z)\\
  &=&\sum_{j=1}^tY(h_ju,z)k_jv+\sum_{i=1}^s\sum_{p,q}Y((a^i_{pq}+b^i_{pq})u,z) e_{pq}^iv\\
   &=&\sum_{p=1}^{n_{i_0}}Y(b_{pq_0}^{i_0}u,z) E_{pq_0}^{i_0}v.
  \end{eqnarray*}
Take $u\in V$ such that $b^{i_0}_{p_0q_0}u\ne0.$  Since  $E_{pq_0}^{i_0}v$ for $p=1,...,n_{i_0}$ are linearly independent, we know from Lemma \ref{2} that $\sum_{p=1}^{n_{i_0}}Y(b_{pq_0}^{i_0}u,z) E_{pq_0}^{i_0}v$ is nonzero. This is a contradiction. The proof is complete. $\hfill\Box$

\begin{Corollary}
  Let $V$ be a simple vertex operator algebra, and $H$  a finite dimensional Hopf algebra. Let $K$ be the kernel of the Hopf action. Then 

  (a) If  $K=0$, then $H$ is a group algebra.

  (b) If $K\neq0$ and $V$ is a semisimple $H$-module, then $K$ is a Hopf ideal and $H/K$ is a group algebra.
  \end{Corollary}

\end{document}